\documentclass[a4paper,11pt]{article}

\newtheorem {theorem}{Theorem}[section]

\newtheorem {lemma}[theorem]{Lemma}

\newcounter{conjecture}
\newcounter{remark}

\newcommand{\eqnsection}{
   \renewcommand{\theequation}{\thesection.\arabic{equation}}
   \makeatletter
   \csname @addtoreset\endcsname{equation}{section}
   \makeatother}

\newtheorem {definition}{Definition}[section]

\author{Alan Hammond\\Department of Statistics\\University of California, Berkeley\\alanmh@stat.berkeley.edu}
\title{Percolation and lattice animals: exponent relations, and conditions for $\theta(p_c)=0$}

\usepackage{amsfonts}
\usepackage[dvips]{graphicx}

\newcommand{\sdash}{\sigma_{n,m}}
\newcommand{\sigm}{\sigma_{n,m}p^n(1-p)^m}
\newcommand{\signt}{\sigma_{n,m}t_n^n(1-t_n)^m}
\newcommand{\sigt}{{\Big( \frac{n}{n+m} \Big)}^{n} {\Big( \frac{m}{n+m} \Big)}^{m}}
\newcommand{\sigo}{p_c^n(1 - p_c)^m}
\newcommand{\reta}{\varsigma}
\newcommand{\rrho}{\varrho}
\newcommand{\appak}{\lambda}
\newcommand{\lapha}{\alpha}

\newcommand{\suf}[2]{\sum_{n}{\sum_{#1}{#2}}}
\newcommand{\suh}[3]{\sum_{n}{\sum_{#1}^{#2}{#3}}}

\newcommand{\atwo}{m > \lfloor n \lapha + n^{1/2} \rfloor + 1}
\newcommand{\alphigh}{m \in n B ( \lapha , C ( \log{n} / n )^{1/2} ) }
\newcommand{\alplow}{m \in n B ( \lapha , n^{-1/2} ) }
\newcommand{\alpnhigh}{m \in n B ( \lapha_n , C ( \log{n} / n )^{1/2} ) }
\newcommand{\alpnlow}{m \in n B ( \lapha_n , n^{-1/2} ) }

\begin{document}
\maketitle

{\bf Abstract}
We examine the percolation model on $\mathbb{Z}^d$ by an approach
involving lattice animals and their surface-area-to-volume ratio.
For $\beta \in [0,2(d-1))$, let
$f(\beta)$ be the asymptotic exponential rate in the number
of edges of the number of lattice animals containing the origin which
have surface-area-to-volume ratio $\beta$. The function $f$ is bounded
above by a function which may be written in an explicit form. For low
values of $\beta$ (\mbox{$\beta \leq 1/p_c - 1$}), equality holds, as
originally demonstrated by F.Delyon. For higher values ($\beta > 1/p_c
- 1$), the inequality is strict. 

We introduce two critical exponents, one of which describes how
quickly $f$ falls away from the explicit form as $\beta$ rises from
$1/p_c - 1$, and the second of which describes how large clusters
appear in the marginally subcritical regime of the percolation model.
We demonstrate that the pair of exponents must satisfy certain
inequalities, while other such inequalities yield sufficient
conditions for the absence of an infinite cluster at the critical
value. The first exponent is related to one of a more conventional nature in
the scaling theory of percolation, that of correlation size. In
deriving this relation, we find that there are two possible
behaviours, depending on the value of the first exponent, for the
typical surface-area-to-volume ratio of an unusually large cluster in
the marginally subcritical regime.

This paper provides an account of the central aspects of the approach,
including the proofs of the main results. In the longer report
\cite{techrep}, complete proofs of all of the assertions are given. 


\section{Introduction}

Percolation on the integer lattice $\mathbb{Z}^d$ is one of the most
fundamental and intensively studied models in the rigorous theory of
statistical mechanics. Many aspects of the behaviour of the model in
the subcritical and supercritical regime have been determined
rigorously. The problem of understanding the behaviour of the model at
criticality, and interplay between this behaviour and that for
parameter values nearby, has been addressed widely by physicists, but
the search for proofs of many of their predictions continues. These
predictions typically take the form of asserting the value of critical
exponents, and thereby describe the power-law decay or explosion of
characteristics of the model near criticality.   

In this paper, we examine the percolation model by an
approach involving lattice animals, divided according to their
surface-area-to-volume ratio.
Throughout, we work with the bond percolation
model in $\mathbb{Z}^d$. However, the results apply to the site or
bond model on any infinite transitive amenable graph with inessential changes.

For any given $p \in (0,1)$, two lattice animals with given size are
equally likely to arise as the cluster $C(0)$ containing the origin
provided that they have the same surface-area-to-volume ratio.
For given $\beta \in (0,\infty)$, there is an exponential growth rate
in the number of edges  for the number of lattice animals up to translation that have surface-area-to-volume ratio very close to $\beta$. This growth rate $f(\beta)$ may be studied as a function of $\beta$.
To illustrate the connection between the percolation model and the combinatorial question of the behaviour of $f$, note that the probability that the cluster containing the origin contains a large number $n$ of edges is given by
\begin{displaymath}
\mathbb{P}_p(\vert C(0) \vert = n) = \sum_{m}{\sigma_{n,m}p^n(1-p)^m},
\end{displaymath}
where $\sigma_{n,m}$ is the number of lattice animals that contain the origin, have $n$ edges and $m$ outlying edges. We rewrite the right-hand-side to highlight the role of the surface-area-to-volume ratio, $m/n$:
\begin{equation}\label{rrr}
\mathbb{P}_p(\vert C(0) \vert = n) = \sum_{m}{(f_n(m/n)p (1-p)^{m/n})^{n}}.
\end{equation}
Here $f_n(\beta) = (\sigma_{n,\lfloor \beta n \rfloor})^{1/n}$ is a
rescaling that anticipates the exponential growth that occurs. We
examine thoroughly the link between percolation and combinatorics
provided by Equation \ref{rrr}. 

An overview of the approach is now given, in the form of a
description of the organisation of the paper.
In Section \ref{secttwo}, we describe the model, and define notations, before
stating the combinatorial results that we will use. The proofs are
largely omitted, as are a few results in later sections. (We refer the
interested reader to the report \cite{techrep}, in which all proofs are given
in full, along with some notes on the literature.) The
combinatorial results assert the existence of the function $f$ and
describe aspects of its behaviour, Theorem \ref{thmtwo} implying that
\begin{equation}\label{fgh}
\log f(\beta) \leq (\beta + 1)\log(\beta + 1) - \beta \log \beta \
\textrm{for $\beta \in (0,2(d-1))$}.
\end{equation}
F.Delyon \cite{DelyonT} showed that equality holds for $\beta \in (0,1/p_c -
1)$. Theorem \ref{thmtwo} implies that the inequality is strict for higher values of
$\beta$. The marked change, as $\beta$ passes through $1/p_c - 1$, in the structure of large
lattice animals of surface-area-to-volume ratio $\beta$ is a
combinatorial analogue of the phase transition in percolation at criticality. 
The notion of a collapse transition for animals has been explored in \cite{MR1304212}.

In Section \ref{sectthr}, two scaling hypotheses are
introduced, each postulating the existence of a critical exponent. One of the
exponents, $\reta$, describes how quickly $f(\beta)$ drops away from
the explicit form given on the right-hand-side of [\ref{fgh}] as $\beta$ rises above $1/p_c - 1$. The other,
$\appak$, describes how rapidly decaying in $n$ is the discrepancy between the
critical value and that value on the subcritical interval at which the
probability of observing an $n$-edged animal as the cluster to
which the origin belongs is maximal. The first main result, Theorem
\ref{thmthr}, is then proved: the inequalities $\appak < 1/2$ and $\reta
\appak < 1$ cannot both be satisfied, because they imply that the mean
cluster size is uniformly bounded on the subcritical interval,
contradicting known results.

In Section \ref{sectfour}, sufficient conditions for the absence of an
infinite cluster at the critical value are proved. Theorem
\ref{thmfour} asserts that $\reta < 2$ or $\appak > 1/2$ are two such
conditions.
Except for some borderline cases, the range of values remaining after
Theorems \ref{thmthr} and \ref{thmfour} is specified by $\appak < 1/2$ and $\reta \appak >
1$. In Theorem \ref{thmund}, where we see that in this case, such a sufficient
condition may be expressed in terms of the extent to which the asymptotic
exponential rate $f(\beta)$ is underestimated by its finite
approximants $f_n(\beta)$ for a certain range of values of $\beta$.
The extent of underestimation is related to combinatorial exponents
such as the entropic exponent (see, for example, \cite{MR95m:82076}).

In Section \ref{sectfive}, we relate the value of $\reta$ to an
exponent of a more conventional nature in the scaling theory of
percolation, that of
correlation size (see Theorem \ref{thmfive}). 
Suppose that we perform an experiment in which the
surface-area-to-volume ratio of the cluster to which the origin
belongs is observed, conditional on its having a very large number of
edges, for a $p$-value slightly below $p_c$. 
How does the typical measurement, $\beta_p$, in this experiment behave 
as $p$ tends to $p_c$? 
The value $\beta_p$ tends to lie somewhere on the interval
\mbox{$(1/p_c - 1,1/p - 1)$}.
In Theorem \ref{thmsix}, we determine that there are two possible scaling
behaviours. The inequality $\reta < 2$ again arises, distinguishing the
two possibilities. If $\reta < 2$, then $\beta_p$ scales much closer to
$1/p_c - 1$ while if $\reta > 2$, it is found
to be closer to $1/p - 1$.

\section{Notations and combinatorial results}\label{secttwo}

Throughout, we work with the bond percolation model on $\mathbb{Z}^d$,
for any given $d \geq 2$.
This model has a parameter $p$ lying in the interval $[0,1]$. Nearest
neighbour edges of $\mathbb{Z}^d$ are declared to be open with
probability $p$, these choices being made independently between
distinct edges. For any vertex $x \in \mathbb{Z}^d$, there is a
cluster $C(x)$ of edges accessible from $x$, namely the collection of
edges that lie in a nearest-neighbour path of open edges one of whose members contains
$x$ as an endpoint. The percolation probability $\theta(p)$ as a
function of $p$ may then be written \mbox{$\theta(p)= \mathbb{P}(
\vert C(0) \vert = \infty)$}. To demonstrate the continuity of
$\theta$, it suffices to show that $\theta(p_c) = 0$ (cf \cite{grim}), where $p_c$
denotes the critical value, namely the infimum of those values of $p$
for which $\theta$ is positive.

\begin{definition}
A lattice animal is the collection of edges of a finite connected
subgraph of $\mathbb{Z}^d$. An edge of $\mathbb{Z}^d$ is said to be
outlying to a lattice animal if it is not a member of the animal, and
if there is an edge in the animal sharing an endpoint with this
edge. 
We adopt the notations:
\begin{itemize}
\item for $n,m \in \mathbb{N}$, set $\Gamma_{n,m}$ equal to the collection
of lattice animals in $\mathbb{Z}^d$ one of whose edges contains the
origin, having $n$ edges, and $m$ outlying edges. Define
$\sigma_{n,m} = \vert \Gamma_{n,m} \vert$. The surface-area-to-volume
ratio of any animal in $\Gamma_{n,m}$ is said to be $m/n$.
\item
for each $n \in \mathbb{N}$, define the function $f_n: [0,\infty) \to [0,\infty)$ by
\begin{displaymath}
 f_n(\beta)=(\sigma_{n,\lfloor \beta n \rfloor})^{1/n}
\end{displaymath} 
\end{itemize}
\end{definition}
On another point of notation, we will sometimes write the index set of
a sum in the form $n S$, with $S \subseteq (0,\infty)$, by which is meant
$\{m \in \mathbb{N}: m/n \in S \}$. \\
We require some results about the asymptotic exponential growth rate
 of the number of lattice animals as a function of their
 surface-area-to-volume ratio. The proofs of the theorems stated here
 are given in \cite{techrep}.
\begin{theorem} \label{thmone} 
${}$
\flushleft
\begin{enumerate}
\item For $\beta \in [0,\infty) - \{ 2(d-1) \}$, $f(\beta)$ exists,
being defined as the limit $\lim_{n \to \infty} f_n(\beta) $. 
\item   for $ \beta > 2(d-1)$, $f(\beta) = 0$. 
\item for $\beta \in (0,2(d-1)), n \in
\mathbb{N}$, $f_n$ satisfies $f_n(\beta) \leq L^{1/n} n^{1/n}
 f(\beta)$, where the constant $L$ may be chosen uniformly in $\beta \in (0,2(d-1))$.
\end{enumerate}
\end{theorem}
\begin{theorem}\label{thmtwo}
${}$
\flushleft
\begin{enumerate}
\item f is log-concave on the interval $(0,2(d-1))$.
\item Introducing $ g: (0,2(d-1)) \to [0,\infty)$ by means of the formula
\begin{displaymath}
 f(\beta) = g(\beta)\frac{ (\beta +1)^{\beta +1} } { \beta^\beta },
\end{displaymath}
we have that
\begin{displaymath}
g(\beta) \left\{ \begin{array}{ll}
        = 1 & \textrm{on $(0,\lapha]$,}\\
        < 1 & \textrm{on $(\lapha,2(d-1))$,} 
\end{array} \right. 
\end{displaymath} 
where throughout $\lapha$ denotes the value $ 1/p_c -1 $.
\end{enumerate}       
\end{theorem}
{\bf Remark} The assertion that $g = 1$ on $(0,\lapha]$ was originally
proved by Delyon \cite{DelyonT}. We include here the proof of the
 other part of the theorem. \\
{\bf Proof}
We must show that, for $\beta \in (\lapha,2(d-1))$, $g(\beta)$ is strictly less than one.
Let $\beta$ lie in this interval. Let $p = 1/(1 + \beta)$. Note that $p < p_c$, and that
\begin{eqnarray}
\mathbb{P}_p(\vert C(0) \vert = n) & \geq & \mathbb{P}_p(C(0) \in \Gamma_{n,\lfloor \beta n \rfloor}) \nonumber \\
 & = & \vert \Gamma_{n,\lfloor \beta n \rfloor} \vert \frac{{\beta}^{\lfloor \beta n \rfloor}}{(1 + \beta)^{n + \lfloor \beta n \rfloor}} \nonumber \\
 & = &  \big( f_n(\beta) \big)^n \frac{{\beta}^{\lfloor \beta n \rfloor}}{(1 + \beta)^{n + \lfloor \beta n \rfloor}}. \nonumber 
\end{eqnarray}
Taking logarithms yields
\begin{displaymath}
 \frac{\log \mathbb{P}_p \big( \vert C(0) \vert = n \big)}{n} \geq
 \log f_n(\beta) + \frac{\lfloor \beta n \rfloor \log \beta}{n} -
 \Big( 1 + \frac{\lfloor \beta n\rfloor }{n} \Big) \log (1 + \beta) ,
\end{displaymath}
from which it follows that
\begin{equation}\label{pop}
 \liminf_{n \to \infty}{\frac{\log \mathbb{P}_p \big( \vert C(0) \vert = n \big)}{n}} \geq  \log f(\beta) + \beta \log \beta - (1+ \beta) \log (1+ \beta) .
\end{equation}
The right-hand-side of (\ref{pop}) is equal to $\log g(\beta)$, by definition.
The exponential decay rate for the probability of observing a large
 cluster in the subcritical phase was established in \cite{MR86h:82045}. Since $p < p_c$, this means the
 left-hand-side of (\ref{pop}) is negative. This implies that $g(\beta)
 < 1$, as required. $\Box$ 

\section{Critical exponents and inequalities}\label{sectthr}

We introduce two scaling hypotheses, each of which proposes the existence of a
critical exponent. We then state and prove the first
main theorem, which demonstrates that a pair of inequalities involving
the two exponents cannot both be satisfied. \\
{\bf Hypothesis $(\appak)$} 
\begin{definition}
For each $n \in \mathbb{N}$, let $t_n \in (0,p_c)$ denote the least
value satisfying the condition
\begin{equation}\label{eteen}
\sum_{m}{\sdash t_n^n (1 - t_n)^m} = \sup_{p \in (0,p_c]}{\sum_{m}{\sdash p^n (1 - p)^m}}. 
\end{equation}
\end{definition}
That is, $t_n$ is some point at or below the critical value at which the probability of observing an $n$-edged animal as the cluster to which the origin belongs is maximal. It is reasonable to suppose that $t_n$ is slightly less than $p_c$, and that the difference decays polynomially in $n$ as $n$ tends to infinity.

\begin{definition} 
Define $\Omega_{+}^{\appak} = \{ \beta \geq 0 : \liminf_{n \to \infty}{(p_c - t_n)/n^{-\beta}} = \infty \}$, and 
$\Omega_{-}^{\appak} = \{ \beta \geq 0 : \limsup_{n \to \infty}{(p_c - t_n)/n^{-\beta}} = 0 \}$. \\
If $\sup{\Omega_{-}^{\appak}} = \inf{\Omega_{+}^{\appak}}$, then hypothesis $(\appak)$ is said to hold, and $\appak$ is defined to be equal to the common value.
\end{definition}
So, if hypothesis $(\appak)$ holds, then $p_c - t_n$ behaves like
$n^{-\appak}$, for large $n$. We remark that it would be consistent
with the notion of a scaling window about criticality that the probability of observing
the cluster $C(0)$ with $n$-edges achieves its maximum on the
subcritical interval on a short plateau whose right-hand endpoint is
the critical value. If this is the case, then $t_n$ should lie at the
left-hand endpoint of the plateau. To be confident that $p_c - t_n$ is
of the same order as the length of this plateau, the definition of the
quantities $t_n$ could be changed, so that a
small and fixed constant multiples the right-hand-side of  
(\ref{eteen}). In this paper, any proof of a statement involving the
exponent $\appak$ is valid if it is defined in terms of this altered
version of the quantities $t_n$. 
 
{\bf Hypothesis $(\reta)$}

This hypothesis is introduced to describe the behaviour of $f$ for
values of the argument just greater than $\lapha$.  Theorem
\ref{thmtwo} asserts that the value $\lapha$ is the greatest for which $\log f(\beta) = (\beta +1)\log(\beta + 1)- \beta \log \beta$; the function $g$ was introduced to describe how $\log f$ falls away from this function as $\beta$ increases from $\lapha$. Thus, we phrase hypothesis $(\reta)$ in terms of $g$.
\begin{definition} 
Define $\Omega_{-}^{\reta} = \{ \beta \geq 0 : \liminf_{\delta \to 0}{(g(\lapha + \delta)- g(\lapha))/{\delta}^{\beta}} = 0 \}$, and 
$\Omega_{+}^{\reta} = \{ \beta \geq 0 : \limsup_{\delta \to
0}{(g(\lapha + \delta) - g(\lapha))/{\delta}^{\beta}} = - \infty \}$. \\
If $\sup{\Omega_{-}^{\reta}} = \inf{\Omega_{+}^{\reta}}$, then hypothesis $(\reta)$ is said to hold, and $\reta$ is defined to be equal to the common value.
\end{definition}
If hypothesis $(\reta)$ holds, then greater values of $\reta$ correspond to a smoother behaviour of $f$ at $\lapha$. For example, if $\reta$ exceeds $N$ for $N \in \mathbb{N}$, then $f$ is $N$-times differentiable at $\lapha$.
\begin{theorem}\label{thmthr}
Suppose that hypotheses ($\reta$) and ($\appak$) hold.  If $\appak < 1/2$, then $\reta \appak \geq 1$.
\end{theorem}
{\bf Proof}
We prove the Theorem by contradiction, assuming that the two
hypotheses hold, and that $\appak < 1/2$, $ \reta \appak < 1$. We will arrive at the conclusion that the mean cluster size, given by
$\sum_{n}{n \mathbb{P}_{p}{(\vert C(0) \vert = n)}}$, is bounded
above, uniformly for $p \in (0,p_c)$. That this is not so is proved in \cite{MR86h:82045}. 
Note that
\begin{displaymath}
\sup_{p \in (0,p_c)}{\sum_{n}{n \mathbb{P}_{p}{(\vert C(0) \vert = n)}}} \leq \sum_{n}{n \mathbb{P}_{t_n}{(\vert C(0) \vert = n)}}.
\end{displaymath}
We write 
\begin{equation}\label{enst}
\mathbb{P}_{t_n}{(\vert C(0) \vert = n)} = \sum_{m}{\sigma_{n,m} t_n^n(1 - t_n)^m},
\end{equation}
and split the sum on the right-hand-side of (\ref{enst}). To do so, we use the following definition.
\begin{definition}\label{alphn}
For $n \in \mathbb{N}$, let $\lapha_n$ be given by $t_n = 1/(1 + \lapha_n)$.
For $G \in \mathbb{N}$, let  $D_n ( = D_n(G))$ denote the interval 
\begin{displaymath}
D_n = ( \lapha_n - G{\{\log(n)/n \}}^{1/2} , \lapha_n + G{\{\log(n)/n \}}^{1/2}).
\end{displaymath} 
\end{definition}
Now, 
\begin{displaymath} 
\sum_{m}{\sigma_{n,m} t_n^n(1 - t_n)^m} = C_1(n) + C_2(n) + C_3(n) ,
\end{displaymath}
where the terms on the right-hand-side are given by 
\begin{eqnarray}
 C_1(n) & = & \sum_{m \in n D_n}{\sigma_{n,m}t_n^n(1 - t_n)^m} , \nonumber \\
 C_2(n) & = &  \sum_{m \in n \big( (0,2(d-1)) - D_n
 \big)}{\sigma_{n,m}t_n^n(1 - t_n)^m} \nonumber
\end{eqnarray}
and
\begin{displaymath}
 C_3(n)  =   \sum_{m \in \{ 2(d-1)n ,\ldots, 2(d-1)n + 2d \}}{\sigma_{n,m}t_n^n(1 - t_n)^m}. 
\end{displaymath}
\begin{definition}\label{defnphi}
Let the function $\phi: (0,\infty)^2 \to \mathbb{R}$ be given by
\begin{displaymath}
\phi (\lapha,\beta) = (\beta + 1)\log(\beta + 1) - \beta \log \beta + \beta \log \lapha - (\beta + 1)\log(\lapha + 1).
\end{displaymath}
\end{definition}
{\bf Remark.} That $\phi \leq 0$ is straigtforward.
\begin{lemma}\label{lemfour}
The function $\phi$ satisfies
\begin{displaymath}
 \phi \big( \lapha , \lapha + \gamma \big) = - \frac{\gamma^2}{2\lapha(\lapha + 1)} + O(\gamma^3).
\end{displaymath}
\end{lemma}
The trivial proof is omitted. \\
We have that
\begin{eqnarray}
 \sum_{n}{C_2(n)}
 & = & \suf{m \in n((\lapha,2(d-1)) - D_n)}{{\bigg( f_n(m/n)\frac{\lapha_n^{m/n}}{(1 + \lapha_n)^{1 + m/n}} \bigg)}^n} \nonumber \\
 & \leq & L \sum_{n}{ n  \sum_{m \in n((\lapha,2(d-1)) - D_n)}{\exp{\{n\phi_{\lapha_n,m/n}\}}}}, \nonumber 
\end{eqnarray}
where the inequality is valid by virtue of Theorem \ref{thmone} and
the fact that $g \leq 1$. Lemma \ref{lemfour} implies that
\begin{displaymath}
 \sum_{m \in n((\lapha,2(d-1)) - D_n)}{\exp{\{n\phi_{\lapha_n,m/n}\}}} \leq (2(d-1) - \lapha) n^{-K} , 
\end{displaymath}
where $K$ may be chosen to be arbitrarily large by an appropriate choice of $G$. It is this consideration that determines the choice of $G$.
The miscellaneous term $C_3$ is treated by the following lemma.
\begin{lemma}\label{lemmanex}
There exists $r \in (0,1)$ such that, for $n$ sufficiently large and for  $m \in \{ 2(d-1)n, \ldots, 2(d-1)n + 2d \}$, we have that
\begin{displaymath} 
\sigma_{n,m} \leq  \frac{ ( 1+ \frac{m}{n} )^{n+ m}}{(\frac{m}{n})^{m}} r^n . 
\end{displaymath}
\end{lemma} 
{\bf Proof} See \cite{techrep}.\\
We find that the $m$-indexed summand in $C_3(n)$ is at most $r^n
\exp{n \phi_{\lapha_n,m/n}}$: thus $C_3(n) \leq (2d+1) r^n$.
Note that $C_1$ satisfies
\begin{eqnarray}
C_1(n) & = & \sum_{m \in n D_n^{*}}{{\bigg( f_n(m/n)\frac{\lapha_n^{m/n}}{(1 + \lapha_n)^{1 + m/n}} \bigg)}^n} \nonumber \\
& \leq &  L n \sum_{m \in n D_n}{g(m/n)^n \exp (n \phi_{\lapha_n,m/n})}, \nonumber 
\end{eqnarray}
where the inequality is a consequence of Theorems \ref{thmone} and \ref{thmtwo}. The fact that the function $\phi$ is nowhere positive implies that
\begin{displaymath}
 C_1(n) \leq L n \sum_{m \in n D_n}{g(m/n)^n}.
\end{displaymath}
Hence the desired contradiction will be reached if we can show that
\begin{equation}\label{expf}
\sum_{n}{ n \sum_{m \in n D_n}{g(m/n)^n}} 
\end{equation}
is finite.
As such, the proof is completed by the following lemma.
\begin{lemma}
Assume hypotheses $(\reta)$ and $(\appak)$. Suppose that $\appak < 1/2$ and that $\reta \appak < 1$. Then, for $\epsilon \in (0,1 - \reta \appak)$ and $n \in \mathbb{N}$ sufficiently large,
\begin{equation}\label{wrte}
\sum_{m \in n D_n}{g(m/n)^n} \leq \exp{-n^{1 - \reta \appak - \epsilon}}.
\end{equation}
\end{lemma}
{\bf Proof}
Let ${\reta}^* > \reta$ and ${\appak}^* > \appak$ be such that
$\appak^* < 1/2$ and ${\reta}^* {\appak}^* < \reta \appak + \epsilon$.
By hypothesis $(\reta)$, there exists $\epsilon'>0$ such that
\begin{displaymath}  
\delta \in (0,{\epsilon}') \ \textrm{implies} \ g(\lapha + \delta) - g(\lapha) < -  {\delta}^{{\reta}^*}.
\end{displaymath}
From Theorems \ref{thmone} and \ref{thmtwo}, it follows that $\sup_{\beta \in [\lapha +
\epsilon',2(d-1)]}{g(\beta)} < 1$, which shows that the contribution
to the sum in (\ref{wrte}) from all those terms indexed by $m$ for
which $m/n > \lapha + \epsilon'$ is exponentially decaying in
$n$. Thus, we may assume that
there exists $N_1$ such that for $n \geq N_1$, if $m \in D_n^{*}$ then $m/n - \lapha < \epsilon'$. Note that, by hypothesis $(\appak)$, $\lapha_n - \lapha \geq n^{-{\appak}^*}$ for sufficiently large. Hence, there exists $N_2$ such that, for $n \geq N_2$,
\begin{displaymath}
\lapha_n - G(\log(n)/n)^{1/2} \geq \lapha + n^{-{\appak}^*} - G(\log(n)/n)^{1/2} \geq \lapha + (1/2)n^{-{\appak}^*}.
\end{displaymath}
For $n \geq \max\{N_1,N_2\}$ and $m \in n D_n^*$,
\begin{eqnarray}
 g(m/n) & \leq & 1 - (m/n - \lapha)^{{\reta}^*} \nonumber \\
        & \leq & 1 - (\lapha_n - G(\log(n)/n)^{1/2} - \lapha)^{{\reta}^*} \nonumber \\
        & \leq & 1 - ((1/2)n^{-{\appak}^*})^{{\reta}^*}. \nonumber
\end{eqnarray}
So, for $n \geq max(N_1,N_2)$,
\begin{displaymath}
\sum_{m \in n D_n^{*}}{g(m/n)^n} \leq (2G(n \log(n))^{1/2})[1 - C' n^{-{\appak}^* {\reta}^*}]^{n},
\end{displaymath}
for some constant $C' > 0$.
There exists $g \in (0,1)$, such that for large $n$, 
\begin{displaymath}
 [1 - C' n^{-{\appak}^* {\reta}^*}]^{n} \leq g^{n^{1 - {\appak}^*{\reta}^*}}.
\end{displaymath}
This implies that
\begin{displaymath}
\sum_{m \in n D_n^{*}}{g(m/n)^n} \leq  h^{n^{1 - {\appak}^*{\reta}^*}} \ \textrm{for large $n$ and $h \in (g,1)$}.
\end{displaymath}
From ${\reta}^* {\appak}^* < \reta \appak + \epsilon$, we find that
\begin{displaymath}
\sum_{m \in n D_n^{*}}{g(m/n)^n} \leq \exp{-n^{1 - \reta \appak - \epsilon}} \ \textrm{for large $n$,}
\end{displaymath}
as required. $\Box$
\section{Sufficient conditions for $\theta(p_c)=0$}\label{sectfour}
In this section, we give two theorems, demonstrating sufficient
conditions for the continuity of the percolation probability in terms
of inequalities on $\reta$ and $\appak$.
\begin{theorem}\label{thmfour}
Assume that hypotheses ($\reta$) and ($\appak$) hold.
\flushleft
\begin{enumerate}
\item Suppose that $\reta < 2$. Then $\theta(p_c) = 0$.
\item Suppose that $\appak > 1/2$. Then $\theta(p_c) = 0$.
\end{enumerate}
\end{theorem}
The proof of Theorem \ref{thmfour} will exploit the characterisation of continuity provided by the following lemma.
\begin{definition}
${}$
\flushleft
\begin{itemize}
\item Let $\sigma(p)= \suf{m}{\sigm}$.
\item Let $\sigma_N(p)=\sum_{n \leq N}{\sum_{m}{\sigm}}$
\end{itemize}
\end{definition}
\begin{lemma}\label{lemghe}
A necessary and sufficient condition for $\theta(p_c)=0$ is that $\sigma_n$ tends uniformly to $\sigma$ on the interval $(0,p_c)$.
\end{lemma}
{\bf Proof} See \cite{techrep}. \\ 
{\bf Proof of Theorem \ref{thmfour}}
By Lemma \ref{lemghe}, to establish that $\theta(p_c)=0$, it suffices
to show that $\sigma_n$ tends to $\sigma$ uniformly on $(0,p_c)$. We
begin by verifying this condition under the hypotheses of the first
part of the Theorem. We will show that 
\begin{equation}\label{escc}
\sum_{n}{\sum_{m}{\sigma_{n,m}\sup_{p \in (0,p_c)}p^n(1-p)^m}} < \infty .
\end{equation}
This will do because
\begin{eqnarray}
\sup_{p \in (0,p_c)}{\big( \sigma(p) - \sigma_{N}(p) \big)} & = & \sup_{p \in (0,p_c)}{\sum_{n \geq N+1}{\sum_{m}\sigm}} {}\nonumber\\
 & \leq & \sum_{n \geq N+1}{\sum_{m}{\sdash \sup_{p \in (0,p_c)}{p^n(1-p)^m}}} {}\nonumber
\end{eqnarray}
So the condition stated in (\ref{escc}) implies the uniform convergence of $\sigma_n$ to $\sigma$ on the subcritical interval.

Note that
\begin{displaymath}
\sup_{p \in (0,p_c)}{p^n(1-p)^m} = \left\{ \begin{array}{ll} \sigt & \textrm{for $n/(n+m) \leq p_c$} \\
\sigo & \textrm{for other pairs $(n,m)$}. \end{array} \right. 
\end{displaymath}

This observation allows us to decompose the sum appearing in (\ref{escc}):

\begin{eqnarray}
\suf{m}{\sdash \sup_{p \in (0,p_c)}{p^n(1-p)^m}} & = & \suh{m  =
1}{\lfloor n \lapha \rfloor }{\sigma_{n,m} \sigo }  \\ \label{thisa}
& + & \suf{m > \lfloor n \lapha \rfloor}{\sigma_{n,m}\sigt} .{}\nonumber
\end{eqnarray}
Now,
\begin{displaymath}
 \suh{m =  1}{\lfloor n \lapha \rfloor}{\sdash \sigo } \leq \suf{m}{\sdash \sigo},
\end{displaymath}
which is less than or equal to one, being the critical probability that the origin lies in a finite cluster.

Set $A$ equal to the second sum on the right-hand-side of (\ref{thisa}).
It suffices to show that $A$ is finite. Our strategy is to split each of the summands of $n$ into two parts, each of which is a sum over $m$ in an interval which has an $n$-dependence. The first sum, $A_1$, will include those $m$-values sufficiently close to $n \lapha$ that this term can be bounded in terms of the critical probability of observing a large cluster. The second sum, $A_2$, will be shown to decay quickly, under the assumption that $\reta < 2$.

Write $A = A_1 + A_2$, where
\begin{eqnarray} 
A_1 & = & \suh{m =  \lfloor n \lapha \rfloor + 1}{\lfloor n \lapha + n^{1/2} \rfloor +1 }{\sdash\sigt}, \nonumber \\
\textrm{and} \, \, A_2 & = & \suf{m > \lfloor n \lapha + n^{1/2} \rfloor +1}{\sdash\sigt}.\nonumber
\end{eqnarray}
Recalling that $\lapha = 1/p_c - 1$,
\begin{displaymath}
A_1 =  \suh{m = \lfloor n \lapha \rfloor + 1}{\lfloor n
\lapha + n^{1/2} \rfloor + 1 }{\sdash\sigo \exp\big( - n \phi (\lapha,m/n) \big)}, 
\end{displaymath}
where the function $\phi$ was specified in Definition \ref{defnphi}.
For each $m \in \{ \lfloor n \lapha \rfloor, \ldots, \lfloor n \lapha
+ n^{1/2} \rfloor + 1 \}$, $c_m \in (0,3/2)$, where $c_m$ is given by $m/n = \lapha + c_m n^{-1/2}$.
Lemma \ref{lemfour} implies that for any sufficiently large $C'$, there exists $N_1$ such that for all $n \geq N_1$, and for $m \in \{ \lfloor n \lapha \rfloor + 1, \ldots, \lfloor n \lapha + n^{1/2} \rfloor + 1 \}$,
\begin{displaymath}
 - \phi (\lapha,m/n) \leq 9/[8 n \lapha(\lapha +1)] + C'/{n^{3/2}} .
\end{displaymath}
From this, we deduce that for $n \geq N_1$ and $m \in \{
\lfloor n \lapha \rfloor + 1, \ldots, \lfloor n \lapha + n^{1/2}
\rfloor + 1 \}$, \mbox{$\exp{(- n \phi (\lapha,m/n) )}$} is bounded above, by $C$, say.
So,
\begin{eqnarray}
A_1 & \leq & \sum_{n < N_1}{\sum_{m \in \{ \lfloor n \lapha \rfloor + 1, \ldots, \lfloor n \lapha + n^{1/2} \rfloor + 1 \}}{\sdash\sigt}} {} \nonumber \\
   {} & & + C   \sum_{n \geq N_1}{\sum_{m \in \{ \lfloor n \lapha \rfloor + 1, \ldots, \lfloor n \lapha + n^{1/2} \rfloor + 1 \}}{\sdash\sigo}}, \nonumber 
\end{eqnarray}
which is finite, as desired.

We now seek to bound $A_2$:
\begin{eqnarray}
A_2 & = & \suf{\atwo}{{\bigg(f_n(m/n)\Big(\frac{n}{n+m}\Big){\Big(\frac{m}{n+m}\Big)}^{m/n}\bigg)}^n} {} \nonumber \\          
    & \leq & L \sum_{n}{\sum_{m = \lfloor n \lapha + n^{1/2} \rfloor +
    2}^{2(d-1)n - 1}{  n
    \bigg(f(m/n) \Big(\frac{n}{n+m}\Big){\Big(\frac{m}{n+m}\Big)}^{m/n}\bigg)^n }} \nonumber \\
   & & + {}  \sum_{n}{\sum_{m = 2(d-1)n}^{2(d-1)n + 2d}{\sdash\sigt}} \nonumber
\end{eqnarray}
where the inequality follows from Theorem \ref{thmone} and the fact
that $g \leq 1$. By Lemma \ref{lemmanex}, there exists $r \in (0,1)$ such that, for $n$ sufficiently large,
\begin{displaymath}
 \sum_{m = 2(d-1)n}^{2(d-1)n + 2d}{\sdash\sigt} \leq (2d+1) r^n . 
\end{displaymath}
It follows from the definition of the function $g$ that
\begin{equation}\label{eat}
A_2 \leq L  \sum_{n}{\sum_{m = \lfloor n \lapha + n^{1/2} \rfloor + 2}^{2(d-1)n - 1}{n g(m/n)^n}} + (2d+1) \sum_{n}{r^n}.
\end{equation}
To bound the first term in the expression on the
right-hand-side of (\ref{eat}),
let $\epsilon \in (0,2 - \reta)$. Let ${\delta}' > 0$, be such that, for $\delta \in (0,{\delta}')$, $g(\lapha + \delta) - g(\lapha) < -{\delta}^{\reta + \epsilon}$. Let $\gamma \in (0,1)$ be such that 
\begin{displaymath}
\sup_{\beta \in (\lapha + {\delta }',2(d-1))}{g(\beta)} < \gamma .
\end{displaymath}
Note that 
\begin{eqnarray}
     &  & \suh{m = \lfloor n \lapha + n^{1/2} \rfloor
     +2}{\lfloor n(\lapha + {\delta}') \rfloor}{ n g(m/n)^n} {} \nonumber \\
 & \leq & \suh{m = \lfloor n \lapha + n^{1/2} \rfloor
 +2}{\lfloor n(\lapha + {\delta}') \rfloor}{ n \Big( 1 - ( m/n -
 \lapha)^{\reta + \epsilon} \Big)^n} {} \nonumber \\
    & \leq & {\delta}' \sum_{n}{n^2 \big( 1 - n^{ - \frac{\reta +
    \epsilon}{2}} \big)^n}. \nonumber 
\end{eqnarray}
Since $\reta + \epsilon < 2$, this expression is finite.
Note also that
\begin{displaymath}
  \suh{m  = \lfloor n(\lapha + {\delta}') \rfloor +
 1}{2(d-1)n - 1}{n g(m/n)^n}  
  \leq  2(d-1)\sum_{n}{n^2 \gamma^n} < \infty . 
\end{displaymath}
We deduce that $A_2$ is finite and in doing so, complete the proof of
the first part of Theorem \ref{thmfour}.

We now prove the second part of the Theorem.
A sufficient condition for continuity is
\begin{equation}\label{anothscc}
\suf{m}{\sdash t_n^n (1 - t_n)^m} < \infty .
\end{equation}
Indeed, the supremum over $p$ in $(0,p_c)$ of $\sigma - \sigma_N$ is
bounded above by the expression in (\ref{anothscc}) with the sum in $n$
being taken over values exceeding $N - 1$. By Lemma \ref{lemghe}, if
(\ref{anothscc}) holds, then $\theta(p_c)=0$.

The fact that $t_n \leq p_c$ implies that $t_n^n(1 - t_n)^m \leq
p_c^n(1 - p_c)^m$ provided that $n/(n+m) > p_c$, which holds if and
only if $m \leq \lfloor n \lapha \rfloor$. From this, we may deduce that
\begin{eqnarray} 
 \suh{m = 1}{\lfloor n \lapha \rfloor}{\sdash t_n^n(1 - t_n)^m} 
& \leq &  \suh{m = 1}{\lfloor n \lapha \rfloor}{\sdash \sigo} \nonumber \\
& \leq &  \suf{m}{\sdash \sigo} \leq 1 \nonumber
\end{eqnarray}
To verify the condition in (\ref{anothscc}), we must bound the expression 
\begin{equation}\label{wret}
\suh{m = \lfloor n \lapha \rfloor + 1}{2(d-1)n + 2d}{\sdash t_n^n(1 - t_n)^m}. 
\end{equation}
To do so, we make the following definition.
\begin{definition}
For $G \in \mathbb{N}$, let  $D_n^{*} ( = D_n^{*}(G))$ denote the interval 
\begin{displaymath}
D_n^{*} = (\max{\{ \lapha, \lapha_n - G{\{\log(n)/n \}}^{1/2} \}}, \lapha_n + G{\{\log(n)/n \}}^{1/2}),
\end{displaymath} 
where the constants $\{ \lapha_n : n \in \mathbb{N} \}$ were specified
in Definition \ref{alphn}.
\end{definition}
Allowing that $G$ will be determined slightly later, we write the
expression in (\ref{wret}) in the form
\begin{eqnarray}
&  & \suf{m \in n D_n^{*}}{\sdash t_n^n(1 - t_n)^m} + \suf{m \in n((\lapha,2(d-1)) - D_n^{*})}{\sdash t_n^n(1 - t_n)^m} \nonumber \\
& & \, + \,  \sum_{n}{\sum_{m = 2(d-1)n}^{2(d-1)n + 2d}{\sdash t_n^n(1 - t_n)^m}} . \label{fvg}
\end{eqnarray}
An argument identical to that by which the term $C_2$ was bounded in the
proof of Theorem \ref{thmthr} yields
\begin{displaymath}
 \suf{m \in n((\lapha,2(d-1)) - D_n^{*})}{\sdash t_n^n(1 - t_n)^m}
 \leq \sum_{n}{n^{-K}}, 
\end{displaymath}
where $K$ may be chosen to be arbitrarily large by an appropriate
choice of $G$, thereby determining how $G$ is chosen. The third term
in (\ref{fvg}) was labelled $C_3(n)$ in the proof of Theorem \ref{thmthr}
and was shown to be bounded above by $(2d+1) r^n$ for $n$ sufficiently
high.
We have that
\begin{eqnarray}
& & \suf{m \in n D_n^{*}}{\sdash t_n^n(1 - t_n)^m} \nonumber \\
& = & \suf{m \in n D_n^{*}}{\sdash \frac{\lapha^m}{(1 + \lapha)^{n+m}}
 \exp{n\Phi(\lapha_n,\lapha,m/n)} }, \label{kzc}
\end{eqnarray}
where 
\begin{eqnarray}
\Phi(\gamma,\lapha,\beta) & = & \beta \log \gamma - (\beta + 1) \log ( \gamma + 1) - \beta \log \lapha + (\beta + 1)\log(\lapha + 1) \nonumber \\
  & = & \beta \log (1 + (\gamma - \lapha)/\lapha ) - (\beta + 1) \log (1 + (\gamma - \lapha)/(1+ \lapha)) \nonumber \\
  & = & - \, \,  \frac{(\gamma - \lapha)^2}{2} [\beta/{\lapha}^2 - (\beta + 1)/(1 + \lapha)^2] \nonumber \\
 & &  + \, \, \frac{(\gamma - \lapha)(\beta - \lapha)}{\lapha(\lapha + 1)} \,  + \, O\big[(\gamma - \lapha)^3 \big]. \nonumber
\end{eqnarray}
We are supposing that hypothesis $(\appak)$ holds, and that $\appak > 1/2$. Let ${\appak}'$ satisfy $\appak > {\appak}' > 1/2$.
In this context,
\begin{displaymath}
\Phi(\lapha_n,\lapha,\beta) =  \frac{(\lapha_n - \lapha)(\beta - \lapha)}{\lapha(\lapha + 1)} - \frac{(\lapha_n - \lapha)^2}{2} [\beta/{\lapha}^2 - (\beta + 1)/(1 + \lapha)^2] + O(n^{-3{\appak}'}).
\end{displaymath}
Now, $\beta \in D_n^{*}$ implies that there exists $C' > 0$ such that
$\beta - \lapha \leq C'n^{-{\appak}'} +C'(\log(n)/n)^{1/2}$ ; since
$\appak' > 1/2$, we may write $\beta - \lapha \leq
C'(\log(n)/n)^{1/2}$, where the value of $C'$ has been increased if
necessary. For such $\beta$, $\Phi(\lapha_n,\lapha,\beta) \leq
C'n^{-{\appak}' - {1/2}}{\log(n)}^{1/2} + n^{-2{\appak}'} +
O(n^{-3{\appak}'})$. This implies that, for all $n$ and $\beta \in
D_n^{*}$, $\exp{ n \Phi(\lapha_n,\lapha,\beta)} < C'$, where once again the value of $C'$ may have changed.
Recalling that $\lapha = 1/p_c - 1$, we deduce from (\ref{kzc}) that
\begin{eqnarray}
 & & \suf{m \in n D_n^{*}}{\sdash t_n^n(1 - t_n)^m} \nonumber \\
 & \leq & C' \suf{m \in n D_n^{*}}{\sdash \sigo} \nonumber \\
 & \leq & C' \suf{m}{\sdash \sigo} \leq C', \nonumber
\end{eqnarray}
proving the second part of Theorem \ref{thmthr}.$\Box$ \\
We now examine the case where
 $\appak < 1/2$ and $\reta \appak > 1$.
\begin{definition}
Let $n \in \mathbb{N}$, and $\beta \in (0,2(d-1))$. Set 
\begin{equation}\label{eanb}
a_n(\beta) = {\bigg( \frac{f_n'(\beta)}{f(\beta)} \bigg) }^n.
\end{equation}
\end{definition}
{\bf Remark} The quantities $a_n(\beta)$ appear in the factorisation of $\sigma_{n,\lfloor \beta n \rfloor}$,
\begin{displaymath}
 \sigma_{n,\lfloor \beta n \rfloor} = a_n(\beta)g(\beta)^n \bigg( \frac{(\beta + 1)^{\beta + 1}}{\beta^\beta} \bigg)^n.
\end{displaymath}
As such, they measure the extent to which the exponential growth rate
$f(\beta)$ is underestimated by $\sigma_{n,m}$.
 
Performing a similar analysis to that undertaken during each part of
Theorem \ref{thmfour} yields the following result. Its proof appears in \cite{techrep}.
\begin{theorem}\label{thmund}
Assume that hypotheses $(\reta)$ and $(\appak)$ hold. Suppose that \mbox{$\appak < 1/2$} and \mbox{$\reta \appak > 1$}. Let $K$ be large. Then there exist constants $\epsilon > 0$ and $C > 0$ such that for each $n \in \mathbb{N}$,
\begin{eqnarray}
 \epsilon \sum_{\alplow}{a_n(m/n)} & \leq & \sum_{m}{\sdash \sigo} \label{edf} \\
                                   & \leq & \sum_{\alphigh}{a_n(m/n)}
				   + \ n^{-K}  \nonumber
\end{eqnarray}
and
\begin{eqnarray}
 \epsilon \sum_{\alpnlow}{a_n(m/n)} & \leq & \sum_{m}{\signt} \label{edg} \\
                                    & \leq  &
				    \sum_{\alpnhigh}{a_n(m/n)} \ + \ n^{-K} . \nonumber 
\end{eqnarray}
\end{theorem}
{\bf Remark} Here, $B(a,b)$ denotes the interval $(a-b,a+b)$. Note
also that it follows from Theorem \ref{thmund} that the 
condition 
\begin{displaymath}
\sum_{\alpnhigh}{a_n(m/n)} < \infty
\end{displaymath}
implies that $\theta(p_c)=0$, without recourse to scaling hypotheses. In examining this condition, bounds on the entropic exponent are revelant (see \cite{MR95m:82076}).  

\section{Scaling law}\label{sectfive}
In this section, we examine the exponential decay rate in $n$ for the
probability of the 
event $\{ C(0)=n \}$ for $p$ slightly less than $p_c$ by our combinatorial
approach. In doing so, we relate the quantity $\reta$ to the exponent
for correlation size, and see how the scaling behaviour for the typical
surface-area-to-volume ratio of unusually large clusters in the
marginally subcritical regime depends on the value of $\reta$.   
\begin{definition}
Let $q:(0,p_c) \to [0,\infty)$ be given by
\begin{displaymath}
q(p) = \lim_{n \to \infty}{\frac{- \log \mathbb{P}_p(\vert C(0) \vert = n)}{n}}.
\end{displaymath}
 Define 
$\Omega_{+}^{\rrho} = \{ \gamma \geq 0 : \liminf_{p \uparrow p_c}  \frac{q(p)}{(p_c - p)^{\gamma}}  = \infty \}$ \ \textrm{and}\\
$\Omega_{-}^{\rrho} = \{ \gamma \geq 0 : \limsup_{p \uparrow p_c}{ \frac{q(p)}{(p_c - p)^{\gamma}}} = 0 \}$.
If $\sup{\Omega_{-}^{\rrho}} = \inf{\Omega_{+}^{\rrho}}$, then hypothesis $(\rrho)$ is said to hold, and $\rrho$ is defined to be equal to the common value.
\end{definition}
{\bf Remark} The existence of $q$ follows from a standard
subadditivity argument. \\
The quantity $\rrho$ might reasonably be called the exponent for
`correlation size'.
\begin{theorem}\label{thmfive}
There exists $\delta' > 0$ and $p_0 \in (0,p_c)$ such that $p \in
(p_0,p_c)$ implies that $q(p)$ is given by 
\begin{displaymath}
\inf_{\beta \in (\lapha,\lapha + \delta')}{- \log g(\beta) - \phi
\big( 1/p - 1,\beta \big)}.
\end{displaymath}
\end{theorem}
The proof, whose details are given in \cite{techrep}, relies on the
fact that the probability that the cluster $C(0)$ has $n$ edges and
$m$ outlying edges in a percolation with parameter $p$ is given by 
$$
a_n(m/n) \exp{n \Big(  \log g(m/n)    + \phi
\big( 1/p - 1,\beta \big)  \Big)},
$$ 
the first term $a_n(m/n)$ having subexponential decay for large
$n$. 

Theorem \ref{thmfive} allows us to deduce a scaling law that relates the combinatorially defined exponent $\reta$ to one which is defined directly from the percolation model.
\begin{theorem}\label{thmsix}
Assume hypothesis $(\reta)$.
\begin{itemize}
\item Suppose that $\reta \in (1,2)$. Then hypothesis $(\rrho)$ holds and
$\rrho = 2$. 
\item Suppose that $\reta \in (2,\infty)$. Then hypothesis $(\rrho)$
holds and $\rrho = \reta$.
\end{itemize}
\end{theorem}
{\bf Proof}
Suppose that $\reta \in (1,2)$. Choose $\epsilon > 0$ so that $1 < \reta
- \epsilon < \reta + \epsilon < 2$. There exists constants $C_1,C_2>0$
such that, for $p \in (p_0,p_c)$ and $\beta \in (\lapha,\lapha +
\delta')$,
\begin{eqnarray}
 (\beta - \lapha)^{\reta + \epsilon} + C_1 \big( \beta - (1/p - 1)
 \big)^2 & \leq & - \log g(\beta) + - \phi \big( 1/p - 1,\beta \big) \label{rtf} \\
 & \leq &  (\beta - \lapha)^{\reta - \epsilon} + C_2 \big( \beta - (1/p - 1)
 \big)^2. \nonumber
\end{eqnarray}
Applying Theorem \ref{thmfive}, we find that
\begin{equation}\label{rtfo}
 (\beta_p - \lapha)^{\reta + \epsilon} + C_1 \big( \beta_p - (1/p - 1)
 \big)^2  \leq q(p), 
\end{equation}
where $\beta_p \in [\lapha,\lapha +
\delta']$ denotes a value at which the infimum in the interval $[\lapha,\lapha +
\delta']$ of the first term in (\ref{rtf}) is  
attained. 
Let $y_p = 1/p - 1 - \lapha$, and let $\sigma_p$ satisfy $\beta_p = \lapha +
y_p^{\sigma_p}$. Then $\beta_p$ and $\sigma_p$ satisfy
\begin{eqnarray}
 (\reta + \epsilon) (\beta_p - \lapha)^{\reta + \epsilon - 1} & = & - 2 C_1
\big( \beta_p - (1/p - 1) \big) \nonumber \\
 (\reta + \epsilon) y_p^{\sigma_p ( \reta + \epsilon - 1)} & = &  2 C_1
\big( y_p - y_p^{\sigma_p} \big) \label{ghb}
\end{eqnarray}
Since $\beta_p \leq 1/p - 1$, $\sigma_p \geq 1$. From this and
(\ref{ghb}) follows $\liminf_{p \uparrow p_c}{\sigma_p} \geq
1/(\reta + \epsilon - 1)$. Applying (\ref{ghb}) again, we deduce that
 $\lim_{p \uparrow p_c}{\sigma_p} =
1/(\reta + \epsilon - 1)$.
Substituting $\sigma_p$ in (\ref{rtfo}) yields
\begin{displaymath}
 y_p^{\sigma_p (\reta + \epsilon)} + C_1 \big( y_p - y_p^{\sigma_p}
 \big)^2 \leq q(p).
\end{displaymath}
The facts that $\lim_{p \uparrow}{\sigma_p} > 1$ and $\lim_{p
\uparrow}{\sigma_p (\reta + \epsilon)} = (\reta + \epsilon)/(\reta +
\epsilon - 1) > 2$ imply that, for a small constant $c$, \mbox{$c (p_c
- p)^2 \leq q(p)$} for values of $p$ just less than $p_c$. A similar
analysis in which $q(p)$ is bounded below by the infimum on the
interval $[\lapha,\lapha + \delta']$ of the third expression in
(\ref{rtf}) implies that for large $C$,  \mbox{$q(p) \leq C (p_c
- p)^2$}, in a similar range of values of $p$. Thus hypothesis $(\rrho)$ holds, and $\rrho = 2$. 

In the case where $\reta > 2$, let $\epsilon > 0$ be such that $\reta >
2 + \epsilon$. Defining $\sigma'_p$ by \mbox{$\beta_p = 1/p - 1 -
y_p^{\sigma'_p}$}, we find that 
\begin{equation}\label{qew}
 (\reta + \epsilon) 
\big( y_p - y_p^{\sigma'_p} \big)^{\reta + \epsilon - 1} =  2 C_1
 y_p^{\sigma'_p}.
\end{equation}
Note that $\beta_p \geq \lapha$ implies that $\sigma'_p \geq 1$. From
(\ref{qew}), it follows that $\liminf_{p \uparrow
p_c}{\sigma'_p} \geq \reta + \epsilon - 1$.  Since $\reta + \epsilon -
1 > 1$, applying (\ref{qew}) again shows that the limit $\lim_{p \uparrow
p_c}{\sigma'_p}$ exists and infact equals $\reta + \epsilon - 1$. Substituting $\sigma'_p$ in (\ref{rtf}) yields
\begin{displaymath}
 \big( y_p - y_p^{\sigma'_p}
 \big)^{\reta+ \epsilon} + C_1 y_p^{2 \sigma'_p} \leq q(p).
\end{displaymath}
The fact that $\liminf_{p \uparrow p_c}{\sigma'_p} > 1$ implies that
\mbox{$c (p_c - p)^{\reta + \epsilon} \leq q(p)$} for values of $p$
just less than $p_c$. Making use of the
inequality $\reta > 2 + \epsilon$ in considering the infimum
of the third term appearing in (\ref{rtf}) yields in this case 
\mbox{$q(p) \leq C (p_c - p)^{\reta - \epsilon}$} for similar values of
$p$. Thus, since $\epsilon$
 may be chosen to be arbitrarily small, we find that, if $\reta > 2$, then hypothesis $(\rrho)$ holds, and
that $\rrho = \reta$. $\Box$ \\
{\bf Acknowledgements} Financial support was provided by a Domus Graduate Scholarship (Competition B) of Merton College, Oxford.
I would like to thank Terry Lyons 
for stimulating and helpful discussions. I thank
John Cardy, Amir Dembo and Mathew Penrose for their helpful comments.   
\bibliography{pbiblio}
\bibliographystyle{plain}
\end{document}